\documentstyle{article}
\begin{document}
\setlength{\baselineskip}{15pt}
\title{Algebraic decoupling of variables for systems of ODEs of 
quasipolynomial form} 
\author{Benito Hern\'{a}ndez--Bermejo$^{\; *}$ \and 
        V\'{\i}ctor Fair\'{e}n$^{\; 1,*}$} 
\date{}

\maketitle

\noindent {\em Departamento de F\'{\i}sica Fundamental, Universidad Nacional 
de Educaci\'{o}n a Distancia. Senda del Rey S/N, 28040 Madrid, Spain.}

\mbox{}

\begin{center} 
{\bf Abstract}
\end{center}

A generalization of the reduction technique for ODEs recently introduced by 
Gao and Liu is given. It is shown that the use of algebraic methods allows 
the extension of the procedure to much more general flows, as well as the 
derivation of simple criteria for the identification of reducible systems. 

\mbox{}

\mbox{}

\mbox{}


\mbox{}

\mbox{}

\mbox{}

\mbox{}

\mbox{}

\mbox{}

\mbox{}

\mbox{}

\mbox{}

\noindent $^1$ To whom all correspondence should be addressed. E-mail: 
vfairen@uned.es

\noindent $^*$ Present address: Departamento de F\'{\i}sica Matem\'{a}tica y 
Fluidos, Universidad Nacional de Educaci\'{o}n a Distancia. Senda del Rey S/N, 
28040 Madrid, Spain.

\pagebreak
\begin{flushleft}
{\bf 1. Introduction}
\end{flushleft}

The problem of finding first integrals and identifying integrability 
conditions of dynamical systems has deserved a continued effort along 
many decades. The relevance and interest of the problem is reflected in the 
number of different procedures developed for these purposes, such as the Lie 
symmetries method \cite{os1}, Carleman embedding \cite{ks1}, Prelle-Singer 
procedure \cite{ps1} or Painlev\'{e} test \cite{se1}, to cite a sample. 
However, none of the presently known methods can account for the problem in 
its full generality, and only partial answers have been developed. 

Among the different analytic tools available, the Quasipolynomial (QP) 
formalism for ODEs has received increasing attention in the last years. This 
interest was initially centered in integrability properties 
\cite{pm1}--\cite{at1} and canonical forms \cite{bv1,bv3}, but 
applications are also starting to reach different fields such as chemical 
kinetics \cite{vb1}, theoretical biochemistry \cite{bv2}, normal forms 
\cite{sb1} and Hamiltonian systems \cite{bv4}. 

In a recent article \cite{gl2}, Gao and Liu have applied a changing variables 
method (CVM from now on) in order to find first integrals of 3D quadratic 
systems of Lotka-Volterra form. This is done by decoupling one of the 
variables of the initial 3D flow, thus reducing the effective dimension of 
the system in one  unit. Analysis of integrability conditions and 
identification of first integrals is thus a much simpler task in the reduced 
2D system. 

It is worth noting that most transformations employed in \cite{gl2} find 
their place in a natural way within the QP formalism. In this work we 
explore the consequences arising from this fact. As Gao and Liu, we shall 
be primarily concerned with the possibility of reducing a flow into a 
two-dimensional one. In addition to the possibility of finding first 
integrals and integrability conditions already mentioned in \cite{gl2}, it 
should be added that knowledge that a system of dimension three or higher can 
be reduced into a two-dimensional one is interesting in itself because it 
excludes the possibility of chaotic behavior ---a problem to which a 
considerable effort has been devoted recently \cite{fh1} in the case of 3D 
systems. 

We shall demonstrate that the CVM can be completely reformulated in terms of 
the QP formalism. This has four major 
implications: 
\begin{description}
   \item{{\em i)\/}} The first is that the procedure can be made more 
systematic and simpler, since all manipulations can be carried out easily in 
terms of matrix algebra. 

   \item{{\em ii)\/}} The second is that the use of the QP formalism allows 
generalizing the scope of the procedure. Generically, the method allows a 
reduction of one unit in the effective dimension of an $n$-dimensional 
system, with arbitrary $n$. This makes the technique particularly interesting 
for the reduction of 3D flows into two-dimensional ones, thus precluding the  
existence of irregular motion. Most examples will accordingly be on standard 
3D systems. However, such reduction into a two-dimensional flow may also be 
possible for some $n$-dimensional systems (an example is given in Subsection 
3.3). Moreover, we shall demonstrate that the procedure is not limited to 
Lotka-Volterra quadratic systems, but is equally valid for flows with 
much more general nonlinearities. 

   \item{{\em iii)\/}} The third is that the use of matrix algebra leads to 
simple criteria for the identification of reducible systems. 

   \item{{\em iv)\/}} The fourth is that some of the CVM transformations are 
just particular cases of wider transformation families that we characterize. 
We shall see in the examples that, in some cases, different members of those 
families are preferable to the CVM choice.
\end{description}

Before describing the reduction technique, it is convenient to give a short 
account of some relevant features of the QP formalism.

\mbox{}

\begin{flushleft}
{\bf 2. Transformations on QP systems}
\end{flushleft}

We shall begin by briefly recalling those basic properties of QP equations 
that will be necessary in what is to come. We refer the reader to the cited 
literature for further details. 

The starting point of the formalism are QP systems of ODE's of the form:
\begin{equation}
   \label{eq:glv}
   \dot{x}_i = x_{i} \left( \lambda _{i} + \sum_{j=1}^{m}A_{ij}\prod_{k=
      1}^{n}x_{k}^{B_{jk}} \right) , \;\:\;\: i = 1 \ldots n , \;\:\; m \geq n
\end{equation}
where $n$ and $m$ are positive integers, and $A$, $B$ and $\lambda$ are 
$n \times m$, $m \times n$ and $n \times 1$ real matrices, respectively. 
In what follows, $n$ will always denote the number of variables of a QP 
system, and $m$ the number of quasimonomials: 
\begin{equation}
   \prod_{k=1}^{n}x_{k}^{B_{jk}} \; , \;\:\;\: j = 1 \ldots m
\end{equation}

We will assume that $m \geq n$ and that $B$ is of maximal rank, i.e. 
rank($B$) $=n$. If $m<n$ and/or rank($B$) is not maximal, then it can be 
shown \cite{bv3} that the system is redundant and can always be reduced to 
the standard situation $m \geq n$ and rank($B$) $=n$, which is our starting 
assumption.

QP equations (\ref{eq:glv}) are form-invariant under quasimonomial 
transformations (QMTs):
\begin{equation}
   x_{i} = \prod_{k=1}^{n} y _{k}^{C_{ik}} , \;\: i=1,\ldots ,n
   \label{bec}
\end{equation}
for any invertible real matrix $C$. After (\ref{bec}), matrices $B$, $A$, and 
$\lambda$ change to 
\begin{equation}
B' = B \cdot C \;\:, \:\;\: A' = C^{-1} \cdot A \;\: , \:\;\: 
\lambda ' = C^{-1} \cdot \lambda \;\: , \:\;\: 
\end{equation}
respectively, but the QP format is preserved. 

Quasimonomial transformations are complemented by the new-time 
transformations (NTTs) of the form \cite{ha1}:
\begin{equation}
    \label{ntt}
    d \tau \; = \; \xi(x_1, \ldots ,x_n) \; dt 
\end{equation}
where $\xi(x_1, \ldots ,x_n)$ is a smooth function. The most important choice 
for $\xi$ in the QP formalism is \cite{gs2}:
\begin{equation}
   \label{xi1}
   \xi(x_1, \ldots ,x_n) = \prod_{i=1}^n x_i^{\, \beta _i}
\end{equation}
where the $\beta _i$ are real constants. With $\xi$ given by (\ref{xi1}), 
transformation (\ref{ntt}) also preserves the QP format. 

Although we will need in certain specific steps of the procedure some 
additional sets of transformations, it is not necessary to elaborate on them 
now. Therefore, we can proceed to describe the reduction method. For this, 
we shall distinguish three cases of increasing complexity. 

\mbox{}

\begin{flushleft}
{\bf 3. Criteria and algorithms for the reduction of systems}
\end{flushleft}

\noindent {\em 3.1 Case I: $\lambda = 0$ and $m = n$}

In this case we shall see that the reduction in one dimension of the system 
is always possible. The set of ODEs takes the form:
\begin{equation}
   \label{case1}
   \dot{x}_i = x_{i} \left( \sum_{j=1}^{n}A_{ij}\prod_{k=
      1}^{n}x_{k}^{B_{jk}} \right) , \;\:\;\: i = 1 \ldots n 
\end{equation}
We look for a QMT such that for the new QP flow we have: 
\begin{equation}
    \label{bcase1}
      B' = B \cdot C = \left( \begin{array}{cccc} 
                     1 & B'_{12} & \ldots & B'_{1n} \\
                     1 & B'_{22} & \ldots & B'_{2n} \\
                     \vdots & \vdots & \mbox{} & \vdots \\
                     1 & B'_{n2} & \ldots & B'_{nn} 
                       \end{array} \right) \;\: ,
\end{equation}
where the columns 2 to $n$ can be chosen arbitrarily (with the obvious  
restriction rank($B'$) $=n$). Note that the CVM choice \cite{gl2} is a 
particular case of (\ref{bcase1}) with $B'_{ij} = \delta _{ij}$, with 
$2 \leq j \leq n$ and symbol $\delta$ standing for Kronecker's delta. Given 
$B'$, we immediately find from (\ref{bcase1}) that $C$ exists and is unique: 
$C$ $=$ $B^{-1} \cdot B'$. Let $y_i$ denote the variables obtained after the 
QMT of matrix $C$. Then we arrive at the following system:
\begin{equation}
   \dot{y}_i = y_1 y_{i} \left( \sum_{j=1}^{n}A'_{ij} \prod_{k=2}^{n}
      y_{k}^{B'_{jk}} \right) , \;\:\;\: i = 1 \ldots n 
\end{equation}
We now rescale the time variable by means of the NTT:
\begin{equation}
   \label{ntt1}
   d \tau = y_1 dt \;\: ,
\end{equation}
where $t$ is the old time variable and $\tau$ the new one. Let us denote from 
now on the derivative of any function $\chi (\tau)$ of a new time $\tau$ as 
$d\chi /d \tau \equiv \hat{\chi}$. Then after (\ref{ntt1}) we are led to:
\begin{equation}
   \label{case13}
   \hat{y}_i = y_{i} \left( \sum_{j=1}^{n}A'_{ij} \prod_{k=2}^{n}
      y_{k}^{B'_{jk}} \right) , \;\:\;\: i = 1 \ldots n 
\end{equation}
Now notice that the only variables appearing in the r.h.s. of equations 2 to 
$n$ of system (\ref{case13}) are precisely $\{ y_2, \ldots, y_n \}$, i.e. 
$y_1$ has been decoupled. Thus the equation for $y_1$ in (\ref{case13}) 
is a quadrature, and system (\ref{case1}) has been reduced to an 
$(n-1)$-dimensional one in the variables $\{ y_2 , \ldots , y_n \}$ and the 
new time $\tau$.

\mbox{}

\noindent {\em Example of Case I: Euler equations for the free rigid body} 

As an example of Case I we can choose Euler's equations for the free rigid 
body \cite{as1} which are given by:
\begin{eqnarray}
     \dot{x}_1 & = & a _1 x_2 x_3  \nonumber    \\
     \dot{x}_2 & = & a _2 x_1 x_3  \label{ej1}  \\
     \dot{x}_3 & = & a _3 x_1 x_2  \nonumber
\end{eqnarray}
The QP matrices of system (\ref{ej1}) are:
\begin{equation}
     B = \left( \begin{array}{ccc}
            -1  &  1  &  1  \\
             1  & -1  &  1  \\
             1  &  1  & -1  
         \end{array} \right) \; , \;\: 
     A = \left( \begin{array}{ccc}
            a_1 &  0  &  0   \\
             0  & a_2 &  0   \\
             0  &  0  & a_3  
         \end{array} \right) \; , \;\: 
     \lambda = \left( \begin{array}{c} 0 \\ 0 \\ 0 \end{array} \right)
\end{equation}
We now apply a QMT of matrix:
\begin{equation}
     C = \left( \begin{array}{ccc}
             1 & 1/2 & 1/2   \\
             1 &  0  & 1/2   \\
             1 & 1/2 & 0  
         \end{array} \right) 
\end{equation}
The resulting QP system is characterized by matrices:
\begin{equation}
\label{ej12}
    B' = \left( \begin{array}{ccc}
             1  &  0  &  0  \\
             1  &  1  &  0  \\
             1  &  0  &  1  
         \end{array} \right) \; , \;\: 
    A' = \left( \begin{array}{ccc}
            -a_1 &    a_2 & a_3    \\
           2 a_1 & -2 a_2 &  0     \\
           2 a_1 &    0   & -2 a_3  
         \end{array} \right) \; , \;\: 
    \lambda ' = \left( \begin{array}{c} 0 \\ 0 \\ 0 \end{array} \right)
\end{equation}
In this case we have chosen $C$ so as to arrive to a matrix $B'$ of the form 
used in \cite{gl2}. Let $\{ y_1, y_2, y_3 \}$ be the variables of the 
transformed equations corresponding to matrices (\ref{ej12}). If we finally 
perform the NTT $d \tau = y_1 dt$ we arrive to:
\begin{eqnarray}
     \hat{y}_1 & = & y_1(-a_1 + a_2y_2 + a_3y_3)  \nonumber     \\
     \hat{y}_2 & = & y_2(2a_1 - 2a_2y_2)          \label{ej1f}  \\
     \hat{y}_3 & = & y_3(2a_1 - 2a_3y_3)          \nonumber 
\end{eqnarray}
This completes the decoupling procedure. In (\ref{ej1f}) the equation of 
$y_1$ has been reduced to a quadrature, while the equations for $y_2$ and 
$y_3$ do not depend on $y_1$. Moreover, the equations for $y_2$ and $y_3$ 
are also decoupled from each other, and can be integrated straightforwardly. 
The integrability of the system is thus made manifest.

\mbox{}

\noindent {\em 3.2 Case II: $\lambda = 0$ and $m > n$}

We now write the system as:
\begin{equation}
   \dot{x}_i = x_{i} \left( \sum_{j=1}^{m}A_{ij}\prod_{k=1}^{n} 
      x_{k}^{B_{jk}} \right) , \;\:\;\: i = 1 \ldots n , \;\:\; m > n
\end{equation}
We again look for a QMT of matrix $C$ such that:
\begin{equation}
    \label{bcase2}
      B' = B \cdot C = \left( \begin{array}{cccc} 
                     1 & B'_{12} & \ldots & B'_{1n} \\
                     1 & B'_{22} & \ldots & B'_{2n} \\
                     \vdots & \vdots & \mbox{} & \vdots \\
                     1 & B'_{m2} & \ldots & B'_{mn} 
                       \end{array} \right) 
\end{equation}
As in Case I, columns 2 to $n$ of $B'$ in (\ref{bcase2}) can be chosen freely 
in such a way that Rank($B'$) $=n$. Now note that matrices $B$ and $B'$ are 
not square but $m \times n$. This implies that a suitable $C$ may or may not 
exist, depending on the form of $B$. Let us denote by $\tilde{B}$ the 
following $m \times (n+1)$ matrix:
\begin{equation}
       \tilde{B} \equiv \left( \begin{array}{ccccc} 
                     B_{11} & B_{12} & \ldots  & B_{1n} &    1   \\
                     B_{21} & B_{22} & \ldots  & B_{2n} &    1   \\
                     \vdots & \vdots & \mbox{} & \vdots & \vdots \\
                     B_{m1} & B_{m2} & \ldots  & B_{mn} &    1 
                        \end{array} \right) 
\end{equation}
It is not difficult to prove that $C$ exists if and only if Rank($\tilde{B}$) 
$=n$, and in this case it is unique. 

If, according to the matrix criterion Rank($\tilde{B}$) $=n$, a suitable $C$ 
exists, then the rest of the procedure is completely similar to Case I: We 
first perform the QMT of matrix $C$. Let $\{ y_1, \ldots , y_n \}$ be the new 
variables of the transformed system. Then we can factor out $y_1$ in each of 
the equations, and eliminate it later by means of the NTT $d \tau = y_1 dt$. 
The result is the decoupling of $y_1$ and the reduction in one unit of the 
effective dimension of the vector field. 

\mbox{}

\noindent {\em Example of Case II: Halphen system}

We now consider the Halphen equations \cite{hs1,ah1} which describe the 
two-monopole system. We shall write them in the form given in \cite{gn1}: 
\begin{eqnarray}
     \dot{x}_1 & = & x_2 x_3 - x_1 x_2 - x_1 x_3  \nonumber    \\
     \do{x}_2 & = & x_1 x_3 - x_1 x_2 - x_2 x_3  \label{ej2}  \\
     \dot{x}_3 & = & x_1 x_2 - x_1 x_3 - x_2 x_3  \nonumber
\end{eqnarray}
In QP terms we have:
\begin{equation}
\label{mej2}
     B = \left( \begin{array}{ccc}
            -1  &  1  &  1  \\
             1  & -1  &  1  \\
             1  &  1  & -1  \\
             1  &  0  &  0  \\
             0  &  1  &  0  \\
             0  &  0  &  1  
         \end{array} \right) \; , \;\: 
     A = \left( \begin{array}{cccccc}
             1 & 0 & 0 &  0  & -1 & -1 \\
             0 & 1 & 0 & -1  &  0 & -1 \\
             0 & 0 & 1 & -1  & -1 &  0 
         \end{array} \right) \; , \;\: 
     \lambda = \mbox{\bf 0}
\end{equation}
It is a simple task to check that rank($\tilde{B}$) $= 3$ when $B$ is given 
by (\ref{mej2}). Consequently, the system can be reduced. For this purpose we 
may choose a QMT of matrix:
\begin{equation}
     C = \left( \begin{array}{ccc}
             1 & 0 & 0 \\
             1 & 1 & 0 \\
             1 & 0 & 1  
         \end{array} \right) 
\end{equation}
After the QMT, the first variable of the transformed system, say $y_1$, can 
be factored out. If we then eliminate it by means of the NTT (\ref{ntt1}) the 
result is:
\begin{eqnarray}
     \hat{y}_1 & = & y_1 (y_2 y_3 - y_2 - y_3)  \nonumber    \\
     \hat{y}_2 & = & -y_2 + y_2^2 - y_2^2 y_3 + y_3   \label{ej22}  \\
     \hat{y}_3 & = & -y_3 + y_3^2 - y_2 y_3^2 + y_2  \nonumber
\end{eqnarray}
As expected, the first variable is reduced to quadrature and the study of 
system (\ref{ej2}) is reduced to a 2D flow corresponding to the equations 
for $y_2$ and $y_3$ in (\ref{ej22}).

\mbox{}

\noindent {\em 3.3 Case III: $\lambda \neq 0$}

We now focus on QP systems of the general form:
\begin{equation}
   \label{case3}
   \dot{x}_i = x_{i} \left( \lambda _{i} + \sum_{j=1}^{m}A_{ij}\prod_{k=
      1}^{n}x_{k}^{B_{jk}} \right) , \;\:\;\: i = 1 \ldots n , \;\:\; m \geq n
\end{equation}
The line of action now consists in reducing the problem to either Case I or 
II, i.e. in suppressing the $\lambda$ terms. For this we first introduce the 
following new variables:
\begin{equation}
\label{tre}
   y_i = e^{- \lambda _i t} x_i \; , \;\:\;\: i = 1 \ldots n
\end{equation}
Transformation (\ref{tre}) is similar to the one employed in \cite{gl2}, and 
was also used in \cite{bs1} in connection with the construction of canonical 
forms for ODEs. We find:
\begin{equation}
   \label{case32}
   \dot{y}_i = y_{i} \left( \sum_{j=1}^{m}A_{ij} e^{\Gamma _j t}
      \prod_{k=1}^{n}y_{k}^{B_{jk}} \right) , \;\:\;\: i = 1 \ldots n \; ,
\end{equation}
where $\Gamma = B \cdot \lambda$. In order to reduce (\ref{case32}) to Cases 
I or II, the following condition is sufficient:
\begin{equation}
\label{lamb}
         \Gamma _1 = \Gamma _2 = \ldots = \Gamma _m \equiv \gamma
\end{equation}
Provided (\ref{lamb}) holds, we can perform the transformation: 
\begin{equation}
        d \tau = e^{\gamma t} dt 
\end{equation}
We finally arrive to the system:
\begin{equation}
   \label{case34}
   \hat{y}_i = y_{i} \left( \sum_{j=1}^{m}A_{ij} 
      \prod_{k=1}^{n}y_{k}^{B_{jk}} \right) , \;\:\;\: i = 1 \ldots n 
\end{equation}
Thus, if condition (\ref{lamb}) is satisfied, equations (\ref{case3}) can be  
reduced to system (\ref{case34}), which corresponds to Case I if $m=n$ and to 
Case II if $m>n$.

This completes our enumeration of criteria and reduction algorithms. We now 
give two examples corresponding to this last situation. 

\mbox{}

\noindent {\em A first example of Case III: Maxwell-Bloch system}

As an example we may consider the Maxwell-Bloch equations for laser systems. 
In the case of periodic boundary conditions, the equations are \cite{lm1}:
\begin{eqnarray}
     \dot{x}_1 & = & -a_1 x_1 + a_2 x_2                    \nonumber    \\
     \dot{x}_2 & = & -a_3 x_2 + a_2 x_1 x_3                \label{ej3}  \\
     \dot{x}_3 & = & -a_4 (x_3 - x_{30}) - 4 a_2 x_1 x_2   \nonumber
\end{eqnarray}
The QP matrices are:
\begin{equation}
     B = \left( \begin{array}{ccc}
            -1  &  1  &  0  \\
             1  & -1  &  1  \\
             1  &  1  & -1  \\
             0  &  0  & -1 
         \end{array} \right) \; , \;\: 
     A = \left( \begin{array}{cccc}
            a_2 &  0  &   0   &      0       \\
             0  & a_2 &   0   &      0       \\
             0  &  0  & -4a_2 &  a_4 x_{30}
         \end{array} \right) \; , \;\: 
     \lambda = \left( \begin{array}{c} 
                     -a_1 \\ -a_3 \\ -a_4 
               \end{array} \right)
\end{equation}
We first compute $\Gamma$:
\begin{equation}
\label{gam3}
     \Gamma = \left( \begin{array}{c} 
                          a_1 - a_3     \\ 
                     - a_1 + a_3 - a_4  \\  
                     - a_1 - a_3 + a_4  \\
                           a_4 
               \end{array} \right)
\end{equation}
Let us look at the compatibility condition (\ref{lamb}) for $\Gamma$ in 
(\ref{gam3}). If we impose $\Gamma _1 = \Gamma _2 = \Gamma _3$ we 
immediately find:
\begin{equation}
\label{e1}
      2 a_1 = a_3 = a_4
\end{equation}
However, it is not possible to simultaneously verify the last requirement 
$\Gamma _4 = \Gamma _i$, for $i=1,2,3$. Then, system (\ref{ej3}) cannot be 
reduced in general. However, we can follow G\"{u}mral and Nutku \cite{gn1} 
and consider the case in which $x_{30}=0$ in equations (\ref{ej3}). The 
resulting system is given by the following QP matrices:
\begin{equation}
\label{m32}
     B = \left( \begin{array}{ccc}
            -1  &  1  &  0  \\
             1  & -1  &  1  \\
             1  &  1  & -1  
         \end{array} \right) \; , \;\: 
     A = \left( \begin{array}{ccc}
            a_2 &  0  &   0   \\
             0  & a_2 &   0   \\
             0  &  0  & -4a_2 
         \end{array} \right) \; , \;\: 
     \lambda = \left( \begin{array}{c} 
                     -a_1 \\ -a_3 \\ -a_4 
               \end{array} \right)
\end{equation}
Now condition (\ref{lamb}) is satisfied iff (\ref{e1}) holds. The 
parameter values that we have obtained in order to verify equation 
(\ref{lamb})
\begin{equation}
      2 a_1 = a_3 = a_4 \;\: , \;\:\;\: x_{30} = 0 \;\: ,
\end{equation}
are precisely those found after some {\em ad hoc\/} transformations by 
G\"{u}mral and Nutku in \cite{gn1} when characterizing the values of the 
parameters for which the Maxwell-Bloch system (\ref{ej3}) is bi-Hamiltonian. 
Notice that such parameter values arise here in a natural way and allow the 
identification of these integrable cases. 

Then, in what follows we will write in (\ref{m32}):
\begin{equation}
       a_1 \equiv \alpha  \;\: , \;\:\;\: a_3 = a_4 \equiv 2 \alpha
\end{equation}
Thus we can perform transformation (\ref{tre}):
\begin{equation}
      y_1 = e^{\alpha t} x_1 \;\: , \;\:\;\:
      y_2 = e^{2 \alpha t} x_2 \;\: , \;\:\;\:
      y_3 = e^{2 \alpha t} x_3 \;\:\;\: ,
\end{equation}
and then the change $d \tau = e^{- \alpha t} dt$. The result is:
\begin{eqnarray}
     \dot{y}_1 & = &  a_2 y_2        \nonumber    \\
     \dot{y}_2 & = &  a_2 y_1 y_3    \label{ej32}  \\
     \dot{y}_3 & = & -4 a_2 y_1 y_2  \nonumber 
\end{eqnarray}
The QP matrices $A'$ and $B'$ of system (\ref{ej32}) coincide, respectively, 
with $A$ and $B$ in (\ref{m32}), while now $\lambda ' = $ {\bf 0}. Since 
we have $m=n=3$ in (\ref{ej32}), equations (\ref{ej3}) have been reduced to 
Case I of the algorithm. As we know from Subsection 3.1, the reduction to a 
2D flow is always possible in this case. 

According to the procedure for Case I, we first apply to system (\ref{ej32}) 
a QMT of matrix:
\begin{equation}
     C = \left( \begin{array}{ccc}
             1 & 1/2 & 1/2 \\
             2 & 1/2 & 1/2 \\
             2 &  1  &  0  
         \end{array} \right) 
\end{equation}
Let $\{ z_1,z_2,z_3 \}$ be the variables of the transformed system. The 
last step is a NTT $d \tau = a_2 z_1 dt$. The final result is:
\pagebreak
\begin{eqnarray}
     \hat{z}_1 & = & z_1 ( z_2 - 1 )                  \nonumber    \\
     \hat{z}_2 & = & 2 z_2 ( 1 - z_2 - 2 z_3 )        \label{ej33}  \\
     \hat{z}_3 & = & 2 z_3 ( 1 + 2 z_3 )              \nonumber
\end{eqnarray}
Then the first variable is decoupled and we obtain a reduced 2D system. Note 
also that the equation for $z_3$ is directly integrable, so the whole system 
is, in fact, reduced to a one-dimensional problem. 

\mbox{}

\noindent {\em An $n$-dimensional example: Riccati projective systems}

We conclude the examples with the Riccati projective equations which have 
recently deserved some attention in different areas, such as selection 
dynamics \cite{se2} or normal forms \cite{sb1}. These systems are given by: 
\begin{equation}
\label{ej4}   
   \dot{x}_i = \lambda _i x_i + x_i \sum_{j=1}^n a_j x_j  
   , \;\:\;\: i = 1 \ldots n 
\end{equation}
We can follow the steps given in Case III and evaluate $\Gamma$. Thus we 
could simplify the system provided condition (\ref{lamb}) is satisfied, i.e. 
$ \lambda _1 $ $=$ $ \ldots $ $=$ $ \lambda _n$. 

We shall not proceed according to this line of action, however. Instead, we 
shall demonstrate that the techniques described above allow solving 
equations (\ref{ej4}) in general. Since $\lambda \neq$ {\bf 0} in (\ref{ej4}) 
we start, as usual, by applying transformation (\ref{tre}):
\begin{equation}
   y_i = e^{- \lambda _i t} x_i \; , \;\:\;\: i = 1 \ldots n
\end{equation}
The result is:
\begin{equation}
\label{ej42}
   \dot{y}_i = y_i \left( \sum_{j=1}^n a_j e^{\lambda _j t} y_j \right) 
   , \;\:\;\: i = 1 \ldots n 
\end{equation}
In general, we cannot factor out the exponentials in equation (\ref{ej42}). 
In other words, relations (\ref{lamb}) will not be usually satisfied. 
However, this is not an unavoidable difficulty in the case of system 
(\ref{ej42}): We can anyhow perform a QMT of the form (\ref{bcase1}) 
described in Subsection 3.1. The best possibility can be easily seen to be: 
\begin{equation}
\label{qmte4}
              C   = \left( \begin{array}{ccccc} 
                        1   &    0   &    0   & \ldots &    0   \\
                        1   &   -1   &    0   & \ldots &    0   \\
                        1   &    0   &   -1   & \ldots &    0   \\
                     \vdots & \vdots & \vdots & \ddots & \vdots \\
                        1   &    0   &    0   & \ldots &   -1 
                       \end{array} \right) 
\end{equation}
Notice that this QMT corresponds to a different choice to the one considered 
in \cite{gl2}. When we perform the QMT of matrix (\ref{qmte4}) on equations 
(\ref{ej42}) the result is: 
\begin{equation}
\label{ej4f1}
   \dot{z}_1 = z_1^2 \left( a_1 e^{\lambda _1 t} + \sum_{j=2}^n 
               \frac{a_j e^{\lambda _j t}}{z_j} \right)
\end{equation}
and
\begin{equation}
\label{ej4f2}
   \dot{z}_i = 0 \;\: , \;\:\;\: i = 2 \ldots n
\end{equation}
The outcome is that, in its final form (\ref{ej4f1})--(\ref{ej4f2}), the 
integrability of system (\ref{ej4}) is made completely explicit ---in fact, 
equations (\ref{ej4f1})--(\ref{ej4f2}) can be integrated trivially. 

\mbox{}

\begin{flushleft}
{\bf 4. Final remarks}
\end{flushleft}

The QP formalism provides the natural operational framework for the changing 
variables method as given in \cite{gl2}: Not only allows its reformulation in 
simpler matrix terms, but also leads naturally to extensions, for example to 
the case of nonquadratic flows. The use of matrix algebra has also made 
possible the derivation of some simple criteria for the identification of 
reducible systems ---criteria which are quite convenient for practical 
purposes.

It is worth insisting that this kind of approach should be especially 
appropriate in the context of 3D sets of ODEs. However, our treatment has 
been completely general in what concerns to the dimension of the system, 
since the possibility of finding higher-dimensional applications cannot be 
excluded, as our last example illustrates. In any case, the final goal has 
always been the reduction into a 2D flow: When this is possible, chaotic 
dynamics is precluded, and further analysis (on parameter space, for 
instance) can be carried out in much simpler terms. 

\mbox{}

\mbox{}

\begin{flushleft}
{\bf Acknowledgements}
\end{flushleft}

This work has been supported by the DGICYT of Spain (grant PB94-0390) and by 
the E.U. (Esprit WG 24490). B. H. acknowledges a doctoral fellowship from 
Comunidad de Madrid. 

\pagebreak


\begin{thebibliography}{99}
   \bibitem{os1} P. J. Olver, Applications of Lie Groups to Differential 
      Equations, Second Ed. (Springer-Verlag, New York, 1993).
   \bibitem{ks1} K. Kowalski and W.-H. Steeb, Nonlinear Dynamical Systems and 
      Carleman Linearization (World Scientific, Singapore, 1991).
   \bibitem{ps1} M. J. Prelle and M. F. Singer, Trans. Amer. Math. Soc. 279 
      (1983) 215.
   \bibitem{se1} W.-H. Steeb and N. Euler, Nonlinear Evolution Equations and 
      Painlev\'{e} Test (World Scientific, Singapore, 1988).
   \bibitem{pm1} M. Peschel and W. Mende, The Predator--Prey Model 
      (Springer-Verlag, Vienna--New York, 1986).
   \bibitem{bs1} L. Brenig, Phys. Lett. A 133 (1988) 378.
   \bibitem{bg1} L. Brenig and A. Goriely, Phys. Rev. A, 40 (1989) 4119.
   \bibitem{gs1} J. L. Gouz\'{e}, Report INRIA 1308 (1990) 1.
   \bibitem{gl1} A. Goriely and L. Brenig, Phys. Lett. A 145 (1990) 245.
   \bibitem{gs2} A. Goriely, J. Math. Phys. 33 (1992) 2728.
   \bibitem{at1} A. Figueiredo, T. M. Rocha Filho and L. Brenig, J. Math. 
      Phys. 39 (1998) 2929.
   \bibitem{bv1} B. Hern\'{a}ndez--Bermejo and V. Fair\'{e}n, Phys. Lett. A 
      206 (1995) 31.
   \bibitem{bv3} B. Hern\'{a}ndez--Bermejo, V. Fair\'{e}n and L. Brenig, 
       J. Phys. A: Math. Gen. 31 (1998) 2415.
   \bibitem{vb1} V. Fair\'{e}n and B. Hern\'{a}ndez--Bermejo, J. Phys. 
      Chem. 100 (1996) 19023.
   \bibitem{bv2} B. Hern\'{a}ndez--Bermejo and V. Fair\'{e}n, Math. Biosci. 
      140 (1997) 1.
   \bibitem{sb1} S. Louies and L. Brenig, Phys. Lett. A 233 (1997) 184.
   \bibitem{bv4} B. Hern\'{a}ndez--Bermejo and V. Fair\'{e}n, Hamiltonian 
      structure and Darboux theorem for families of generalized 
      Lotka-Volterra systems, J. Math. Phys. (1998), to be published.
   \bibitem{gl2} P. Gao and Z. Liu, Phys. Lett. A 244 (1998) 49.
   \bibitem{fh1} Z. Fu and J. Heidel, Nonlinearity 10 (1997) 1289.
   \bibitem{ha1} J. Hietarinta, B. Grammaticos, B. Dorizzi and A. Ramani, 
      Phys. Rev. Lett. 53 (1984) 1707.
   \bibitem{as1} V. I. Arnold, Mathematical Methods of Classical Mechanics, 
      Second Ed. (Springer-Verlag, New York, 1989).
   \bibitem{hs1} M. Halphen, C. R. Acad. Sci. Paris 92 (1881) 1101.
   \bibitem{ah1} M. Atiyah and N. Hitchin, The Geometry And Dynamics of 
      Magnetic Monopoles (Princeton University Press, Princeton, New Jersey, 
      1988).
   \bibitem{gn1} H. G\"{u}mral and Y. Nutku, J. Math. Phys. 34 (1993) 5691.
   \bibitem{lm1} F. T. Arecchi, in: Order and Chaos in Nonlinear Physical 
      Systems, eds. S. Lundqvist, N. H. March and M. P. Tosi (Plenum Press, 
      New York, 1988) p. 193.
   \bibitem{se2} A. V. Shapovalov and E. V. Evdokimov, Physica D 112 (1998) 
      441.
\end{thebibliography}
\end{document}